\documentclass[12pt, a4paper, hidelinks]{extarticle}

\usepackage[left=30mm, top=20mm, right=15mm, bottom=20mm]{geometry}
\usepackage{setspace}

\usepackage{algorithm}
\usepackage{algpseudocode}
\usepackage{amsmath,amssymb,amsthm}
\usepackage{mathtools}
\usepackage{mathrsfs}
\usepackage{blkarray}
\usepackage{systeme}
\usepackage{indentfirst}
\usepackage{enumitem}
\usepackage[colorlinks=true,citecolor=blue,urlcolor=blue,linkcolor=blue]{hyperref}
\usepackage{cmap}

\usepackage[T2A]{fontenc}
\usepackage[english]{babel}

\usepackage[style=numeric, sorting=nty]{biblatex}
\floatname{algorithm}{Algorithm}

\newtheorem{definition}{Definition}
\newtheorem{theorem}{Theorem}
\newtheorem{proposition}{Proposition}

\newtheorem{corollary}{Corollary}
\newtheorem{remark}{Remark}

\allowdisplaybreaks

\providecommand{\keywords}[1]{
	\small	
	\textbf{\textit{Keywords---}} #1
	\normalsize
}

\title{Construction of a set of circulant matrix submatrices for faster MDS matrix verification}
\author{Stanislav S. Malakhov$ ^1 $}
\date{$ ^1 $National Research University Higher School of Economics\\
	\hfill \break
}

\begin{document}
	
	\maketitle
	\onehalfspacing

	\begin{abstract}
		This paper is obsolete and superseded by the article
		\begin{itemize}
			\item[] Malakhov, S.S. On the circulant matrix MDS testing and the search for circulant MDS matrices. \textit{Cryptogr. Commun.} \textbf{17}, 87–119 (2025). \item[] \href{https://doi.org/10.1007/s12095-024-00746-7}{DOI: /10.1007/s12095-024-00746-7}
		\end{itemize}
	\end{abstract}

	\vspace{6em}
\section{Introduction}
	Suppose a $ k \times m $ matrix $ \boldsymbol{M} $ over a finite field $ \mathbb{F}_{q} $. Then a set
	\begin{displaymath}
		\left\lbrace \left( \boldsymbol{x},\boldsymbol{x} \cdot \boldsymbol{M} \right) | \boldsymbol{x} \in \left( \mathbb{F}_{q} \right)^k \right\rbrace 
	\end{displaymath}	
	is a linear $ \left[ n,k,d \right] $ code of the length $ n = k + m $ and the dimension $ k $ with the minimum Hamming distance $ d $ between any two code words. For a linear $ \left[ n,k,d \right] $ code the Singleton bound holds:
	\begin{displaymath}
		d \leq n - k + 1 = m + 1,
	\end{displaymath}
	and a code with $ d = m + 1 $ is called the maximum distance separable (MDS) code. As revealed in \cite{macwilliams1977theory} on page 321, a linear code is MDS if and only if every square submatrix of $ \boldsymbol{M} $ is non-singular. With respect to an MDS code we shall refer to $ \boldsymbol{M} $ as the MDS matrix.
	\begin{definition}
		A matrix $ \boldsymbol{M} $ is the MDS matrix if every square submatrix of $ \boldsymbol{M} $ is non-singular.
	\end{definition}
	
	MDS matrices are demanded for block cryptographic algorithms where they are responsible for a diffusion quality. The MDS matrix performs a linear transformation of an input block $ \boldsymbol{x} $ of the following property: if $ i, 1 \leq i \leq k, $ elements of $ \boldsymbol{x} $ are altered, then at least $ m - i + 1 $ elements of the output block $ \boldsymbol{x} \cdot \boldsymbol{M} $ alter, where both the input and the output blocks can be interpreted as vectors of a \textit{k}-dimensional vector space over a finite field $ \mathbb{F}_{q} $. In this sense MDS matrices provide perfect diffusion \cite{junod2004}. Several algorithms utilize MDS matrices among which are block ciphers AES, GOST R 34.12---2015, and IDEA NXT and hash functions GOST R 34.11---2012 and Whirlpool.
	
	The two primary objectives regarding MDS matrices are their construction and verification that a particular matrix is MDS one. In general case, both the objectives are computational hard, especially for larger matrices, but there is room for improvement when examining special matrices. % This paper focuses on circulant matrices and proposes some techniques which allow reduction of computational efforts into either circulant MDS matrix construction and verification. 
	This paper focuses on circulant matrices and proposes some techniques which allow reduction of computational efforts into circulant MDS matrix verification. 
	
	\section{Equivalent submatrices of circulant matrices}
	
	Consider circulant matrices $ \boldsymbol{L} $ and $ \boldsymbol{R} $,
	\begin{displaymath}
		\boldsymbol{L} = \begin{pmatrix}
			a_{0} & a_{1} & a_{2} & \cdots &  & a_{n-1} \\
			a_{1} & a_{2} & a_{3} &\cdots & a_{n-1} & a_{0} \\
			a_{2} & a_{3} & \cdots & a_{n-1} & a_{0} & a_{1} \\
			\cdots &  &  &  &  &  \\
			a_{n-2} & a_{n-1} & a_{0} & \cdots &  & a_{n-3} \\
			a_{n-1} & a_{0} & \cdots &  & a_{n-3} & a_{n-2} \\ 		
		\end{pmatrix}, 	\boldsymbol{R} = \begin{pmatrix}
			a_{0} & a_{1} & a_{2} & \cdots & a_{n-2} & a_{n-1} \\
			a_{n-1} & a_{0} & a_{1} & \cdots & a_{n-3} & a_{n-2} \\
			a_{n-2} & a_{n-1} & a_{0} & \cdots &  & a_{n-3} \\
			\cdots &  &  & \ddots &  & \\
			a_{2} & \cdots &  & a_{n-1} & a_{0} & a_{1} \\
			a_{1} & a_{2} & \cdots &  & a_{n-1} & a_{0} \\ 		
		\end{pmatrix},
	\end{displaymath}
	where entries $ a_{0}, \dots, a_{n - 1} $ are pairwise distinct variables.		
	\begin{remark}
		Obviously, one can obtain $ \boldsymbol{L} $ from $ \boldsymbol{R} $ reordering the rows of $ \boldsymbol{R} $ and vice versa. Therefore, it is enough to examine one matrix only, say $ \boldsymbol{L} $. However, for the sake of completeness we shall observe both the matrices where it appears to be appropriate. 
	\end{remark}
	\begin{remark} \label{On elements equality in L and R}
		For $ k \in \{0, \dots, n - 1 \} $ the entries $ a_{(i + k)_{\bmod{n}}(j - k)_{\bmod{n}}} $ of $ \boldsymbol{L} $ equal each other, and so do the entries $ a_{(i + k)_{\bmod{n}}(j + k)_{\bmod{n}}} $ of matrix .
	\end{remark}
	
	Previously, circulant matrices were studied in several works, for instance, in \cite{cauchois2019circulant,gupta2014constructions,kesarwani2019exhaustive,li2016construction}.
	
	\begin{definition} \label{Matrices Equivalence Definition}
		Let $ \boldsymbol{A}_{1} $ and $ \boldsymbol{A}_{2} $ be the matrices of variables. We shall call $ \boldsymbol{A}_{1} $ and $ \boldsymbol{A}_{2} $ equivalent ($ \boldsymbol{A}_{1} \sim \boldsymbol{A}_{2} $), if one can be transformed into the other by reordering the rows and columns.
	\end{definition}

	\begin{remark}
		It is clear that the relation ($ \sim $) is an equivalence relation on the set $ \mathscr{A} $ of matrices, since for each matrix $ \boldsymbol{A}_{1}, \boldsymbol{A}_{2} $ and $ \boldsymbol{A}_{3} $ in $ \mathscr{A} $ the following hold:
		\begin{enumerate}[label = {(\roman*)}]
			\item $ \boldsymbol{A}_{1} \sim \boldsymbol{A}_{1} $,
			\item $ \left\lbrace \boldsymbol{A}_{1} \sim \boldsymbol{A}_{2} \right\rbrace \Longleftrightarrow \left\lbrace \boldsymbol{A}_{2} \sim \boldsymbol{A}_{1} \right\rbrace $,
			\item $ \left\lbrace \boldsymbol{A}_{1} \sim \boldsymbol{A}_{2} \right\rbrace \wedge \left\lbrace \boldsymbol{A}_{2} \sim \boldsymbol{A}_{3} \right\rbrace \Rightarrow \left\lbrace \boldsymbol{A}_{1} \sim \boldsymbol{A}_{3} \right\rbrace $.
		\end{enumerate}
	\end{remark}

	\begin{remark}
		Equivalent matrices may be of interest due to equality of absolute values of their determinants: if $ \boldsymbol{A}_{1} \sim \boldsymbol{A}_{2} $, then $ \left| \det{\boldsymbol{A}_{1}} \right|  = \left| \det{\boldsymbol{A}_{2}} \right| $ whenever exist.
	\end{remark}

	\begin{definition} \label{Submatrix Notation}
		For a matrix $ \boldsymbol{M} $ of a dimension $ n_{1} \times n_{2} $ by $ \boldsymbol{M}\left[ i_{0},\dots,i_{u - 1}; j_{0},\dots,j_{v - 1} \right] $ we denote a submatrix formed from the rows $ i_{0} < \dots < i_{u - 1} $ and columns $ j_{0} < \dots < j_{v - 1} $ of $ \boldsymbol{M} $. Here $ u \in \{1,\dots,n_{1}\} $, and $ v \in \{1,\dots,n_{2}\} $.
	\end{definition}
	
	\begin{definition}
		Let independently $ \boldsymbol{A}_{L} $ and $ \boldsymbol{A}_{R} $ be the respective submatrices of $ \boldsymbol{L} $ and $ \boldsymbol{R} $: 
		\begin{align*}
			\boldsymbol{A}_{L} &= \boldsymbol{L}\left[ i_{0},\dots,i_{u - 1}; j_{0},\dots,j_{v - 1} \right],  \\
			\boldsymbol{A}_{R} &= \boldsymbol{R}\left[ i_{0},\dots,i_{u - 1}; j_{0},\dots,j_{v - 1} \right].
		\end{align*}
		For some $ k \in \{0,\dots,n - 1\} $ by $ \boldsymbol{A}^{(k)}_{L} $ we denote a submatrix of $ \boldsymbol{L} $ formed from $ k $ circular shifts applied in $ \boldsymbol{L} $ to the rows and columns of $ \boldsymbol{A}_{L} $ down and to the left respectively:
		\begin{multline*}
			\boldsymbol{A}^{(k)}_{L} \stackrel{def}{=} \boldsymbol{L}\left[ \pi\left(\left( i_{0} + k \right)_{\bmod{n}}\right) ,\dots,\pi\left(\left( i_{u - 1} + k \right)_{\bmod{n}}\right)\right. ; \\ \left. \rho\left(\left( j_{0} - k \right)_{\bmod{n}}\right) ,\dots, \rho\left(\left( j_{v - 1} - k \right)_{\bmod{n}}\right) \right].
		\end{multline*}
		Here $ \pi $ is a permutations of $ \{(i_{0} + k)_{\bmod n}, \dots,  (i_{u - 1} + k)_{\bmod n}\} $ and $ \rho $ is a permutation of $ \{\left( j_{0} - k \right)_{\bmod{n}}, \dots, \left( j_{v - 1} - k \right)_{\bmod{n}}\} $ such that 
		\begin{align*}
			\pi\left( (i_{0} + k)_{\bmod n}\right) &< \cdots < \pi\left( (i_{u - 1} + k)_{\bmod n}\right), \\
			\rho\left(\left( j_{0} - k \right)_{\bmod{n}}\right) &< \dots < \rho\left(\left( j_{v - 1} - k \right)_{\bmod{n}}\right).
		\end{align*}
		Similarly, for some $ k \in \{0,\dots,n - 1\} $ by $ \boldsymbol{A}^{(k)}_{R} $ we denote a submatrix of $ \boldsymbol{R} $ formed from $ k $ circular shifts applied in $ \boldsymbol{R} $ to the rows and columns of $ \boldsymbol{A}_{R} $ down and to the right respectively:
		\begin{multline*}
			\boldsymbol{A}^{(k)}_{R} \stackrel{def}{=} \boldsymbol{R}\left[ \pi\left(\left( i_{0} + k \right)_{\bmod{n}}\right) ,\dots,\pi\left(\left( i_{u - 1} + k \right)_{\bmod{n}}\right)\right. ; \\ \left. \rho\left(\left( j_{0} + k \right)_{\bmod{n}}\right) ,\dots, \rho\left(\left( j_{v - 1} + k \right)_{\bmod{n}}\right) \right].
		\end{multline*}
		Here $ \pi $ is a permutations of $ \{(i_{0} + k)_{\bmod n}, \dots,  (i_{u - 1} + k)_{\bmod n}\} $ and $ \rho $ is a permutation of $ \{\left( j_{0} + k \right)_{\bmod{n}}, \dots, \left( j_{v - 1} + k \right)_{\bmod{n}}\} $ such that 
		\begin{align*}
			\pi\left( (i_{0} + k)_{\bmod n}\right) &< \cdots < \pi\left( (i_{u - 1} + k)_{\bmod n}\right), \\
			\rho\left(\left( j_{0} + k \right)_{\bmod{n}}\right) &< \dots < \rho\left(\left( j_{v - 1} + k \right)_{\bmod{n}}\right).
		\end{align*}
		In these cases we shall refer to a submatrix $ \boldsymbol{A}_{L} = \boldsymbol{A}^{(0)}_{L} $ or $ \boldsymbol{A}_{R} = \boldsymbol{A}^{(0)}_{R} $ as the generator submatrix, and we shall refer to a submatrix $ \boldsymbol{A}^{(k)}_{L} $ or $ \boldsymbol{A}^{(k)}_{R} $ as the derived submatrix, whenever $ k \in \{1,\dots,n - 1\} $.
	\end{definition}
	
	\begin{proposition} \label{Equivalence Classes}
		Let independently $ \boldsymbol{A}_{L} $ and $ \boldsymbol{A}_{R} $ are respective generator submatrices of the matrix $ \boldsymbol{L} $ and $ \boldsymbol{R} $. Then $ \boldsymbol{A}^{(k)}_{L} \sim \boldsymbol{A}_{L} $ and $ \boldsymbol{A}^{(k)}_{R} \sim \boldsymbol{A}_{R} $ for each $ k \in \left\lbrace 1, \dots, n - 1 \right\rbrace $.
	\end{proposition}
	\begin{proof}
		Consider the case of the matrix $ \boldsymbol{L} $. Let
		\begin{equation*} \label{A_L}
			\boldsymbol{A}_{L} = \begin{pmatrix}
				a_{i_{0}j_{0}} & \cdots & a_{i_{0}j_{s-1}} & a_{i_{0}j_{s}} & \cdots & a_{i_{0}j_{v-1}} \\
				\vdots &  &  &  & \\
				a_{i_{r-1}j_{0}} & \cdots & a_{i_{r-1}j_{s}} & a_{i_{r-1}j_{s}} & \cdots & a_{i_{r-1}j_{v-1}} \\
				a_{i_{r}j_{0}} & \cdots & a_{i_{r}j_{s}} & a_{i_{r}j_{s}} & \cdots & a_{i_{r}j_{v-1}} \\
				\vdots &  &  &  & \\
				a_{i_{u-1}j_{0}} & \cdots & a_{i_{u-1}j_{s}} & a_{i_{u-1}j_{s}} & \cdots & a_{i_{u-1}j_{v-1}} \\
			\end{pmatrix},
		\end{equation*}
		and $ (i_{l} + k)_{\bmod{n}} \leq i_{l} $ for $ l \in \left\lbrace r, \dots, u - 1 \right\rbrace, 0 \leq r \leq u - 1 $, while $ (j_{m} - k)_{\bmod{n}} \leq j_{m} $ for $ m \in \{s, \dots, v - 1\}, 0 \leq s \leq v - 1 $. Then $ \boldsymbol{A}^{(k)}_{L} $ takes the following form
		\begin{multline*} \label{A_L^(k)}
			\left( \begin{matrix}
				a_{(i_{r} + k)_{\bmod{n}}(j_{s} - k)} & \cdots & a_{(i_{r} + k)_{\bmod{n}}(j_{v-1} - k)} \\
				\vdots &  &  \\
				a_{(i_{u-1} + k)_{\bmod{n}}(j_{s} - k)} & \cdots & a_{(i_{u-1} + k)_{\bmod{n}}(j_{v-1} - k)} \\
				a_{(i_{0} + k)(j_{s} - k)} & \cdots & a_{(i_{0} + k)(j_{v-1} - k)} \\
				\vdots &  &  \\
				a_{(i_{r-1} + k)(j_{s} - k)} & \cdots & a_{(i_{r-1} + k)(j_{v-1} - k)} \\
			\end{matrix} \right. \\
			\left. \begin{matrix}
				a_{(i_{r} + k)_{\bmod{n}}(j_{0} - k)_{\bmod{n}}} & \cdots & a_{(i_{r} + k)_{\bmod{n}}(j_{s-1} - k)_{\bmod{n}}} \\
				\vdots &  & \\
				a_{(i_{u-1} + k)_{\bmod{n}}(j_{0} - k)_{\bmod{n}}} & \cdots & a_{(i_{u-1} + k)_{\bmod{n}}(j_{s-1} - k)_{\bmod{n}}} \\
				a_{(i_{0} + k)(j_{0} - k)_{\bmod{n}}} & \cdots & a_{(i_{0} + k)(j_{s-1} - k)_{\bmod{n}}} \\
				\vdots &  & \\
				a_{(i_{r-1} + k)(j_{0} - k)_{\bmod{n}}} & \cdots & a_{(i_{r-1} + k)(j_{s-1} - k)_{\bmod{n}}} \\
			\end{matrix} \right).
		\end{multline*}
%		\begin{displaymath}
%			\begin{pmatrix}
%				a_{(i_{r} + k)_{\bmod{n}}(j_{s} - k)} & \cdots & a_{(i_{r} + k)_{\bmod{n}}(j_{u} - k)} & a_{(i_{r} + k)_{\bmod{n}}(j_{1} - k)_{\bmod{n}}} & \cdots & a_{(i_{r} + k)_{\bmod{n}}(j_{s-1} - k)_{\bmod{n}}} \\
%				\vdots &  &  &  &  & \\
%				a_{(i_{u} + k)_{\bmod{n}}(j_{s} - k)} & \cdots & a_{(i_{u} + k)_{\bmod{n}}(j_{u} - k)} & a_{(i_{u} + k)_{\bmod{n}}(j_{1} - k)_{\bmod{n}}} & \cdots & a_{(i_{u} + k)_{\bmod{n}}(j_{s-1} - k)_{\bmod{n}}} \\
%				a_{(i_{1} + k)(j_{s} - k)} & \cdots & a_{(i_{1} + k)(j_{u} - k)} & a_{(i_{1} + k)(j_{1} - k)_{\bmod{n}}} & \cdots & a_{(i_{1} + k)(j_{s-1} - k)_{\bmod{n}}} \\
%				\vdots &  &  &  &  & \\
%				a_{(i_{r-1} + k)(j_{s} - k)} & \cdots & a_{(i_{r-1} + k)(j_{u} - k)} & a_{(i_{r-1} + k)(j_{1} - k)_{\bmod{n}}} & \cdots & a_{(i_{r-1} + k)(j_{s-1} - k)_{\bmod{n}}} \\
%			\end{pmatrix}.
%		\end{displaymath}
		Note, that $ a_{(i + k)_{\bmod{n}}(j - k)_{\bmod{n}}} = a_{ij} $ in $ \boldsymbol{L} $ for $ k \in \{ 0, \dots, n - 1 \} $, therefore
		\begin{displaymath}
			\boldsymbol{A}^{(k)}_{L} = \begin{pmatrix}
				a_{i_{r}j_{s}} & \cdots & a_{i_{r}j_{v-1}} & a_{i_{r}j_{0}} & \cdots & a_{i_{r}j_{s-1}} \\
				\vdots &  &  &  &  & \\
				a_{i_{u-1}j_{s}} & \cdots & a_{i_{u-1}j_{v-1}} & a_{i_{u-1}j_{0}} & \cdots & a_{i_{u-1}j_{s-1}} \\
				a_{i_{0}j_{s}} & \cdots & a_{i_{0} j_{v-1}} & a_{i_{0}j_{0}} & \cdots & a_{i_{0}j_{s-1}} \\
				\vdots &  &  &  &  & \\
				a_{i_{r-1}j_{s}} & \cdots & a_{i_{r-1}j_{v-1}} & a_{i_{r-1}j_{0}} & \cdots & a_{i_{r-1}j_{s-1}} \\
			\end{pmatrix}.
		\end{displaymath}		
		Since $ \boldsymbol{A}^{(k)}_{L} $ differs from $ \boldsymbol{A}_{L} $ in the order of rows and columns, $ \boldsymbol{A}^{(k)}_{L} \sim \boldsymbol{A}_{L} $. 
		
		The case of the matrix $ \boldsymbol{R} $ may be proved in the similar manner.
	\end{proof}

	\begin{theorem} \label{Equivalence necessary and sufficient condition}
		Suppose $ \boldsymbol{A}_{L} $ and $ \boldsymbol{A}'_{L} $ are submatrices of the matrix $ \boldsymbol{L} $ and independently $ \boldsymbol{A}_{R} $ and $ \boldsymbol{A}'_{R} $ are submatrices of the matrix $ \boldsymbol{R} $. The following statements hold:
		\begin{enumerate}[label = {(\roman*)}]
			\item \label{Theorem 1 Case 1} $ \boldsymbol{A}'_{L} \sim \boldsymbol{A}_{L} $, if and only if $ \boldsymbol{A}'_{L} $ is a derived submatrix from $ \boldsymbol{A}_{L} $;
			\item \label{Theorem 1 Case 2} $ \boldsymbol{A}'_{R} \sim \boldsymbol{A}_{R} $, if and only if $ \boldsymbol{A}'_{R} $ is a derived submatrix from $ \boldsymbol{A}_{R} $.
		\end{enumerate}
	\end{theorem}
	\begin{proof}
		Consider the case \ref{Theorem 1 Case 1}. Necessity immediately follows from proposition \ref{Equivalence Classes}: if $ \boldsymbol{A}'_{L} $ is a derived submatrix from a generator submatrix $ \boldsymbol{A}_{L} $ then $ \boldsymbol{A}'_{L} \sim \boldsymbol{A}_{L} $.
		 
		To prove sufficiency suppose 
		\begin{displaymath}
			\boldsymbol{A}_{L} = \boldsymbol{L}\left[ i_{0},\dots,i_{u - 1}; j_{0},\dots,j_{v - 1} \right],
		\end{displaymath}  
		and $ (i_{l} + k)_{\bmod{n}} \leq i_{l} $ for $ l \in \left\lbrace r, \dots, u - 1 \right\rbrace, 1 \leq r \leq u - 1 $, while $ (j_{m} - k)_{\bmod{n}} \leq j_{m} $ for $ m \in \{s, \dots, v - 1\}, 1 \leq s \leq v - 1 $, so that  
		\begin{multline*}
			\boldsymbol{A}^{(k)}_{L} = \boldsymbol{L}\left[ i_{r} + k - n,\dots,i_{u - 1} + k - n, i_{0} + k, \dots, i_{r - 1} + k;\right. \\ \left. j_{s} - k, \dots, j_{v - 1} - k, n + j_{0} - k, \dots, n + j_{s - 1} - k \right].
		\end{multline*}
		Note, that for a submatrix formed from rows $ i_{0}, \dots, i_{u - 1} $ and columns $ j_{0}, \dots, j_{v - 1} $ of $ \boldsymbol{L} $ there exist $ u $ permutations of rows and $ v $ permutations of columns which arrange another submatrix of $ \boldsymbol{L} $. Indeed, $ r $ may take one value out of $ u $, while $ s $ may take one value out of $ v $ although both $ r $ and $ s $ depend on $ k $. Besides, for $ k \in \{ 1, \dots, n \} $ elements $ a_{(i + k)_{\bmod{n}}(j - k)_{\bmod{n}}} $ of $ \boldsymbol{L} $ are equal, therefore, shifting rows $ i_{0}, \dots, i_{u - 1} $ down and columns $ j_{0}, \dots, j_{v - 1} $ to the left $ k $ times in a circular manner, one can construct all $ u $ permutations of rows and all $ v $ permutations of columns of $ \boldsymbol{A}_{L} $. However, some tuples of rows and columns may repeat. Thereby, all $ u $ permutations of rows and all $ v $ permutations of columns of $ \boldsymbol{A}_{L} $ form the a of submatrices derived from $ \boldsymbol{A}_{L} $. Hence, assuming $ \boldsymbol{A}'_{L} \neq \boldsymbol{A}^{(k)}_{L} $ for each $ k \in \{1,\dots,n - 1\} $, we conclude that neither $ \boldsymbol{A}'_{L} \nsim \boldsymbol{A}_{L} $ nor $ \boldsymbol{A}_{L}' $ is a submatrix of $ \boldsymbol{L} $. The latter case contradicts the assumption that $ \boldsymbol{A}_{L} $ is a submatrix of $ \boldsymbol{L} $, whereupon if $ \boldsymbol{A}'_{L} \neq \boldsymbol{A}_{L} $, then $ \boldsymbol{A}'_{L} \nsim \boldsymbol{A}_{L} $. Thus, $ \boldsymbol{A}'_{L} \sim \boldsymbol{A}_{L} $ if and only if $ \boldsymbol{A}'_{L} $ is a submatrix derived from $ \boldsymbol{A}_{L} $. Now the case \ref{Theorem 1 Case 2} can be proved similarly.
	\end{proof}
	
	It follows from theorem \ref{Equivalence necessary and sufficient condition} that generator matrices $ \boldsymbol{A}_{L} $ and $ \boldsymbol{A}_{R} $ define the respective equivalence classes 
	\begin{align*}
		\left[ \boldsymbol{A}_{L} \right] = \left\lbrace \boldsymbol{A}^{(k)}_{L} | k \in \{ 0, \dots, n - 1 \} \right\rbrace, \\ \left[ \boldsymbol{A}_{R} \right] = \left\lbrace \boldsymbol{A}^{(k)}_{R} | k \in \{ 0, \dots, n - 1 \} \right\rbrace.
	\end{align*} 
	
	\begin{corollary} \label{Corollary of equivalence necessary and sufficient condition}
		Let $ \boldsymbol{A}_{L} $ and $ \boldsymbol{A}'_{L} $ be equivalent submatrices of the matrix $ \boldsymbol{L} $,
		\begin{align*}
			\boldsymbol{A}_{L} &= \boldsymbol{L}\left[ i_{0},\dots,i_{u - 1}; j_{0},\dots,j_{v - 1} \right],  \\
			\boldsymbol{A}'_{L} &= \boldsymbol{L}\left[ i'_{0},\dots,i'_{u - 1}; j'_{0},\dots,j'_{v - 1} \right].
		\end{align*}
		Let also $ \boldsymbol{d} $ and $ \boldsymbol{d}' $ be the tuples, corresponding respectively to $ \boldsymbol{A}_{L} $ and $ \boldsymbol{A}'_{L} $ the following way: The elements of $ \boldsymbol{d} $ and $ \boldsymbol{d}' $ represent differences between ordinal numbers in $ \boldsymbol{L} $ of two consecutive rows of $ \boldsymbol{A}_{L} $ and $ \boldsymbol{A}'_{L} $, i.e. 
		\begin{align*}
			\boldsymbol{d} &= \left( (i_{1} - i_{0}), \dots, (i_{u-1} - i_{u-2}), n - (i_{u-1} - i_{0}) \right),  \\
			\boldsymbol{d}' &= \left( (i'_{1} - i'_{0}), \dots, (i'_{u-1} - i'_{u-2}), n - (i'_{u-1} - i'_{0}) \right).
		\end{align*}
		Then the tuple $ d' $ can be formed by applying circular shifts to $ d $. In other words,
		\begin{displaymath}
			\left( \boldsymbol{A}'_{L} \sim \boldsymbol{A}_{L} \right) \Longrightarrow \left( \exists r \in \{1, \dots, u - 1\}: \boldsymbol{d}' = \sigma^{r-1}(\boldsymbol{d}) \right),
		\end{displaymath}
		where $ \sigma^{r-1}(\boldsymbol{d}) $ denotes a circular shift of the tuple $ \boldsymbol{d} $ to $ r - 1 $ positions to the left, $ r \in \{1, \dots, u - 1\} $.
	\end{corollary}
	\begin{proof}
		It follows from theorem \ref{Equivalence necessary and sufficient condition} that  $ \boldsymbol{A}'_{L} \sim \boldsymbol{A}_{L} $, if and only if $ \boldsymbol{A}'_{L} $ is a submatrix derived from $ \boldsymbol{A}_{L} $. Hence, there exists $ k \in \{0, \dots, n - 1\} $ such that $ \boldsymbol{A}'_{L} = \boldsymbol{A}^{(k)}_{L} $. Thus far,
		\begin{multline*}
			\left( \boldsymbol{A}'_{L} = \boldsymbol{A}^{(k)}_{L} \right) \Longleftrightarrow \left( \boldsymbol{A}'_{L} = \boldsymbol{L}\left[ i_{r} + k - n,\dots,i_{u - 1} + k - n, i_{0} + k, \dots, i_{r - 1} + k;\right.\right.  \\ \left.\left. j_{s} - k, \dots, j_{v - 1} - k, n + j_{0} - k, \dots, n + j_{s - 1} - k \right] \right).
		\end{multline*}
		Consequently,
		\begin{displaymath}
			\begin{split}
				\boldsymbol{d}' =& \left( i_{r + 1} + k - n - (i_{r} + k - n), \dots, i_{u - 1} + k - n - (i_{u - 2} + k - n),\right. \\
				& \quad \left.  i_{0} + k - (i_{u - 1} + k - n), i_{1} + k - (i_{0} + k), \dots, i_{r - 1} + k - (i_{r - 2} + k) \right) \\
				=& \left( (i_{r + 1} - i_{r}), \dots, (i_{u - 1} - i_{u - 2}), n - (i_{u - 1} - i_{0}), (i_{1} - i_{0}), \dots, (i_{r - 1} - i_{r - 2}) \right) \\
				=& \sigma^{r-1}(\boldsymbol{d}).
			\end{split} 
		\end{displaymath}		
	\end{proof}
	
	\begin{remark} \label{Corollary of equivalence necessary and sufficient condition for R}
		By analogy to corollary \ref{Corollary of equivalence necessary and sufficient condition} for equivalent submatrices $ \boldsymbol{A}'_{R} $ and $ \boldsymbol{A}_{R} $ the corresponding tuples $ \boldsymbol{d}' $ and $ \boldsymbol{d} $ can be formed from each other by applying circular shifts.
	\end{remark}
	
	As a result, if $ \boldsymbol{d} $ and $ \boldsymbol{d}' $ cannot be formed from each other by applying circular shifts, then the corresponding submatrices of $ \boldsymbol{L} $ or $ \boldsymbol{R} $ are non-equivalent.
	
%	Оценим возможные мощности классов эквивалентности $ \left[ \boldsymbol{A}_{L} \right] $ и $ \left[ \boldsymbol{A}_{R} \right] $.
	
%	\begin{proposition}
%		Пусть $ \boldsymbol{A}_{L} = \boldsymbol{L}\left[ i_{0},\dots,i_{u - 1}; j_{0},\dots,j_{v - 1} \right] $ --- подматрица матрицы $ \boldsymbol{L} $,	где $ u \leq n \geq v $, $ n $ --- порядок матрицы $ \boldsymbol{L} $. Равенство $ \boldsymbol{A}_{L} = \boldsymbol{A}^{(k)}_{L} $ выполняется, если и только 
%		\begin{enumerate}[label = {(\roman*)}]
%			\item \label{Proposition 2 Case 1} найдется такое $ m \in \{ 1, \dots, u \} $, что $ (i_{(l+m)_{\bmod{u}}} - i_{l})_{\bmod{n}} = k_{1} $ для каждого $ l \in \{0, \dots, u - 1\} $;
%			\item \label{Proposition 2 Case 2} найдется такое $ m' \in \{ 1, \dots, v \} $, что $ (j_{(l'+m')_{\bmod{v}}} - i_{l'})_{\bmod{n}} = k_{2} $ для каждого $ l' \in \{0, \dots, v - 1\} $;
%			\item \label{Proposition 2 Case 3} $ k $ кратно $ k_{1} $ и $ k_{2} $.
%		\end{enumerate}
%	\end{proposition}
	
	\begin{proposition} \label{k-values for A_L}  An $ u \times v $ generator submatrix $ \boldsymbol{A}_{L} $ of matrix $ \boldsymbol{L} $ generates an equivalence class $ \left[ \boldsymbol{A}_{L} \right] $, that contains $ k $ elements,
		\begin{displaymath}
			k \in \left\lbrace qn / \gcd(n, u, v) | q \in \left\lbrace 1, \dots, \gcd(n, u, v) \right\rbrace \right\rbrace.
		\end{displaymath}
	\end{proposition}
	\begin{proof}
		Let $ \boldsymbol{A}_{L} = \boldsymbol{L}\left[ i_{0},\dots,i_{u - 1}; j_{0},\dots,j_{v - 1} \right], $ and $ (i_{l} + k)_{\bmod{n}} \leq i_{l} $ for $ l \in \left\lbrace r, \dots, u - 1 \right\rbrace $,  $ 0 \leq r \leq u - 1 $, while $ (j_{m} - k)_{\bmod{n}} \leq j_{m} $ for $ m \in \{s, \dots, v - 1\}, 0 \leq s \leq v - 1 $. Then the \textit{k}-th submatrix derived 
		\begin{multline*}
		 	\boldsymbol{A}^{(k)}_{L} = \boldsymbol{L}\left[ i_{r} + k - n, \dots, i_{u - 1} + k - n, i_{0} + k, \dots, i_{r - 1} + k;\right. \\ \left.  j_{s} - k, \dots, j_{v - 1} - k, j_{0} + n - k, \dots, j_{s - 1} + n - k \right].
		\end{multline*} 
		
		An equivalence class $ \left[ \boldsymbol{A}_{L} \right] $ contains $ k $ elements, if $ \boldsymbol{A}_{L} = \boldsymbol{A}_{L}^{(k)} $ and $ \boldsymbol{A}_{L} \neq \boldsymbol{A}_{L}^{(k')} $ for each $ k' \in \{1, \dots, k - 1 \} $. Since $ \boldsymbol{A}_{L} $ equals $ \boldsymbol{A}_{L}^{(k)} $ the following systems of equations hold.
		\begin{displaymath}
			\left\lbrace \begin{aligned}
				i_{0} &= i_{r} + k - n; \\ \vdots \\ i_{u - 1 - r} &= i_{u - 1} + k - n; \\ i_{u - r} &= i_{0} + k; \\ \vdots \\ i_{u - 1} &= i_{r - 1} + k,
			\end{aligned} \right. \qquad \left\lbrace \begin{aligned}
				j_{0} &= j_{s} - k; \\ \vdots \\ j_{v - 1 - s} &= j_{v - 1} - k; \\ j_{v - s} &= j_{0} + n - k; \\ \vdots \\ j_{v - 1} &= j_{s - 1} + n - k,
			\end{aligned} \right. 		
		\end{displaymath}
		in variables $ i_{0}, \dots, i_{u - 1}, j_{0}, \dots, j_{v - 1} $ and $ k $. The sums of equations within each system render the next system of Diophantine equations in unknown $ k, r, s $:
		\begin{displaymath}
		\left\lbrace \begin{aligned}
				uk + nr &= nu; \\ vk - ns &= 0,
			\end{aligned} \right. \Longleftrightarrow \left\lbrace \begin{aligned}
				k &= n - q'n / \gcd(n, u); \\ r &= q'u / \gcd(n, u); \\ s &= v - q'v / \gcd(n, u),
			\end{aligned} \right. \Longleftrightarrow \left\lbrace \begin{aligned}
				k &= n - q''n / \gcd(n, v); \\ r &= q''u / \gcd(n, v); \\ s &= v - q''v / \gcd(n, v),
			\end{aligned} \right.
		\end{displaymath}
		where $ q' \in \left\lbrace 1, \dots, \gcd(n, u) \right\rbrace, q'' \in \left\lbrace 1, \dots, \gcd(n, v) \right\rbrace $. Note, that both the forms represent the same solution, producing a new Diophantine equation
		\begin{align*}
			q' / \gcd(n, u) = q'' / \gcd(n, v) &\Longleftrightarrow q' \gcd(n, v) - q'' \gcd(n, u) = 0 \\ &\Longleftrightarrow \left\lbrace \begin{aligned}
				q' &= \gcd(n, u) - \frac{q\gcd(n, u)}{\gcd(\gcd(n, u),\gcd(n, v))}; \\ q'' &= \gcd(n, v) - \frac{q\gcd(n, v)}{\gcd(\gcd(n, u),\gcd(n, v)),}
			\end{aligned} \right. \\ &\Longleftrightarrow \left\lbrace \begin{aligned}
			q' &= \gcd(n, u) - q\gcd(n, u) / \gcd(n, u, v); \\ q'' &= \gcd(n, v) - q\gcd(n, v) / \gcd(n, u, v),	
			\end{aligned} \right.
		\end{align*}
		where $ q \in \left\lbrace 1, \dots, \gcd(n, u, v) \right\rbrace $. Then finally
		\begin{displaymath}
			\left\lbrace \begin{aligned}
				uk + nr &= nu; \\ vk - ns &= 0,	
			\end{aligned} \right. \Longleftrightarrow \left\lbrace \begin{aligned}
				k &= qn / \gcd(n, u, v); \\ r &= u - qu / \gcd(n, u, v); \\ s &= qv / \gcd(n, u, v),
			\end{aligned} \right. 
		\end{displaymath}
		with $ q \in \left\lbrace 1, \dots, \gcd(n, u, v) \right\rbrace $.
	\end{proof}
	
	\begin{remark} \label{Remark on A_R to k-values for A_L}
		For a submatrix $ \boldsymbol{A}_{R} = \boldsymbol{R}\left[ i_{0},\dots,i_{u - 1}; j_{0},\dots,j_{v - 1} \right] $ such that $ (i_{l} + k)_{\bmod{n}} \leq i_{l} $ and $ (j_{m} + k)_{\bmod{n}} \leq j_{m} $, when $ l \in \left\lbrace r, \dots, u - 1 \right\rbrace,  0 \leq r \leq u - 1 $ and $ m \in \left\lbrace s, \dots, v - 1 \right\rbrace,  0 \leq s \leq v - 1 $, the next system of equations holds
		\begin{displaymath}
			\left\lbrace \begin{aligned}
				uk + nr &= nu; \\ vk + ns &= nu,	
			\end{aligned} \right. \Longleftrightarrow \left\lbrace \begin{aligned}
				k &= qn / \gcd(n, u, v); \\ r &= u - qu / \gcd(n, u, v); \\ s &= v - qv / \gcd(n, u, v),
			\end{aligned} \right. 
		\end{displaymath}
		where $ q \in \left\lbrace 0, \dots, \gcd(n, u, v) \right\rbrace $.
	\end{remark}
	
	\begin{remark}
		If there exist $ r $ and $ s $, and the $ u \times u $ submatrix $ \boldsymbol{A}_{L} $, for which $ \boldsymbol{A}_{L} = \boldsymbol{A}^{(k)}_{L} $, then $ \boldsymbol{A}_{L} = \boldsymbol{A}^{(\hat{k})}_{L} $ for $ \hat{k} = n - k, \hat{r} = s $ and $ \hat{s} = r $. Indeed, since $ \hat{r} + \hat{s} = u = r + s $, $ \hat{s} = r $ for $ \hat{r} = s $, and also
		\begin{displaymath}
			\left\lbrace \begin{aligned}
				\hat{k} &= \hat{q}n / \gcd(n,u); \\ \hat{r} &= u - \hat{q}u / \gcd(n, u); \\ \hat{r} &= qu / \gcd(n, u),
			\end{aligned} \right. \Rightarrow \left\lbrace \begin{aligned}
				\hat{k} &= \hat{q}n / \gcd(n,u); \\ \hat{q} &= \gcd(n,u) - q, 
			\end{aligned} \right. \Rightarrow \hat{k} = n - nq / \gcd(n,u) = n - k. 
		\end{displaymath}
		
		In particular, this means that either $ k = n $ or $ k \leq n / 2 $, while $ q \in \{0, \dots, \gcd(n, u) / 2 \} $.
	\end{remark}
	
	\begin{remark} \label{Remark on prime order matrix submatrices quantity}	
		It is obvious that for prime $ n $ every submatrix of $ \boldsymbol{L} $ or $ \boldsymbol{R} $ generates an equivalence class that contains exactly $ n $ elements, unless the dimension of the submatrix is $ n \times n $.
	\end{remark}

\section{Transposes and anti-transposes of submatrices of circulant matrices}
	
	Not only the equivalent matrices, conforming definition \ref{Matrices Equivalence Definition}, have equal absolute values of their determinants but also the transposed ones. The following proposition formulates the necessary and sufficient condition at which one submatrix equals the transpose of another submatrix.
	
	\begin{proposition} \label{On Matrix Transpose}
		Let $ \boldsymbol{A} = \boldsymbol{L}\left[ i_{0},\dots,i_{u - 1}; j_{0},\dots,j_{v - 1} \right] $ be a submatrix of the matrix $ \boldsymbol{L} $ of the dimension $ n \times n $, where $ u \leq n \geq v $. Then $ \boldsymbol{A}^T = \boldsymbol{L}\left[ j_{0},\dots,j_{v - 1}; i_{0},\dots,i_{u - 1} \right] $.
	\end{proposition}
	\begin{proof}
		\begin{displaymath}
			\boldsymbol{A} = \begin{pmatrix}
				a_{i_{0}j_{0}} & \cdots & a_{i_{0}j_{s}} & \cdots & a_{i_{0}j_{v - 1}} \\
				\vdots & & & & \\
				a_{i_{r}j_{0}} & \cdots & a_{i_{r}j_{s}} & \cdots & a_{i_{r}j_{v - 1}} \\ 
				\vdots & & & & \\
				a_{i_{u - 1}j_{0}} & \cdots & a_{i_{u - 1}j_{s}} & \cdots & a_{i_{u - 1}j_{v - 1}} \\ 
			\end{pmatrix} \Rightarrow \boldsymbol{A}^T \stackrel{def}{=} \begin{pmatrix}
				a_{i_{0}j_{0}} & \cdots & a_{i_{r}j_{0}} & \cdots & a_{i_{u - 1}j_{0}} \\
				\vdots & & & & \\
				a_{i_{0}j_{s}} & \cdots & a_{i_{r}j_{s}} & \cdots & a_{i_{u - 1}j_{s}} \\ 
				\vdots & & & & \\
				a_{i_{0}j_{v - 1}} & \cdots & a_{i_{r}j_{v - 1}} & \cdots & a_{i_{u - 1}j_{v - 1}} \\ 
			\end{pmatrix}
		\end{displaymath}
		Note, that $ a_{ij} = a_{(i + j - 1)_{\bmod{n}}} = a_{ji} $ for each element $ a_{ij} $ of $ \boldsymbol{L} $, therefore
		\begin{align*}
			\boldsymbol{A}^T \stackrel{def}{=}& \begin{pmatrix}
				a_{i_{0}j_{0}} & \cdots & a_{i_{r}j_{0}} & \cdots & a_{i_{u - 1}j_{0}} \\
				\vdots & & & & \\
				a_{i_{0}j_{s}} & \cdots & a_{i_{r}j_{s}} & \cdots & a_{i_{u - 1}j_{s}} \\ 
				\vdots & & & & \\
				a_{i_{0}j_{v - 1}} & \cdots & a_{i_{r}j_{v - 1}} & \cdots & a_{i_{u - 1}j_{v - 1}} \\ 
			\end{pmatrix} = \begin{pmatrix}
				a_{j_{0}i_{0}} & \cdots & a_{j_{0}i_{r}} & \cdots & a_{j_{0}i_{u - 1}} \\
				\vdots & & & & \\
				a_{j_{s}i_{0}} & \cdots & a_{j_{s}i_{r}} & \cdots & a_{j_{s}i_{u - 1}} \\ 
				\vdots & & & & \\
				a_{j_{v - 1}i_{0}} & \cdots & a_{j_{v - 1}i_{r}} & \cdots & a_{j_{v - 1}i_{u - 1}} \\ 
			\end{pmatrix} \\
			=& \boldsymbol{L}\left[ j_{0},\dots,j_{v - 1}; i_{0},\dots,i_{u - 1} \right].
		\end{align*}
		Furthermore, $ \boldsymbol{A}^T $ indeed is a submatrix of $ \boldsymbol{L} $.
	\end{proof}
	\begin{corollary}
		Under the conditions of theorem \ref{On Matrix Transpose}
		\begin{displaymath}
			\left\lbrace \boldsymbol{A} = \boldsymbol{A}^T \right\rbrace \Longleftrightarrow \left\lbrace \{ i_{0},\dots,i_{u - 1} \} = \{ j_{0},\dots,j_{v - 1} \} \right\rbrace. 
		\end{displaymath}
	\end{corollary}
	\begin{proof}
		To prove necessity suppose that $ \boldsymbol{A} = \boldsymbol{A}^T $. In this case $ u = v $, and $ a_{i_{r}j_{s}} = (\boldsymbol{A})_{rs} = \left( \boldsymbol{A}^T\right) _{rs} = a_{j_{r}i_{s}} $, for $ r \in \{0,\dots,u - 1\} $ and $ s \in \{0,\dots,u - 1\} $. As a corollary to remark \ref{On elements equality in L and R} the equality $ a_{i_{r}j_{s}} = a_{j_{r}i_{s}} $ of elements of $ \boldsymbol{L} $ leads to the system of equations
		\begin{displaymath}
			\left\lbrace \begin{aligned}
				i_{r} &= j_{r} + k \pmod n; \\ j_{s} &= i_{s} - k \pmod n,
			\end{aligned} \right. \Rightarrow \left\lbrace \begin{aligned}
				i_{r} &= j_{r}; \\ k &= 0,
			\end{aligned} \right.			
		\end{displaymath}
		Hence, $ \{ i_{0},\dots,i_{u - 1} \} = \{ j_{0},\dots,j_{v - 1} \} $.
		
		To prove sufficiency suppose that $ \{ i_{0},\dots,i_{u - 1} \} = \{ j_{0},\dots,j_{v - 1} \} $. In this case
		\begin{displaymath}
			\boldsymbol{A} = \boldsymbol{L}\left[ i_{0},\dots,i_{u - 1}; j_{0},\dots,j_{v - 1} \right] = \boldsymbol{L}\left[ j_{0},\dots,j_{v - 1}; i_{0},\dots,i_{u - 1} \right] = \boldsymbol{A}^T.
		\end{displaymath}
	\end{proof}

	A statement similar to one making by proposition \ref{On Matrix Transpose} holds also for the matrix $ \boldsymbol{R} $. However, in this case instead of the transpose, one should observe the reflection of the elements along the diagonal $ i + j = v - 1 $, where $ i \in \{0,\dots,u - 1\}, j \in \{0,\dots,v - 1\} $ and $ v $ is the number of columns of a submatrix. 
	
	\begin{definition}
		Given the arbitrary matrix $ \boldsymbol{M} $ of dimension $ n_{1} \times n_{2} $,
		\begin{displaymath}
			\boldsymbol{M} = \left( m_{ij} \right) | i \in \{0,\dots,n_{1} - 1\}, v \in \{0,\dots,n_{2} - 1\},
		\end{displaymath}
		by $ \boldsymbol{M}^\tau $ we denote the matrix of dimension $ n_{2} \times n_{1} $ such that
		\begin{displaymath}
			\boldsymbol{M}^\tau = \left( m_{(n_{1} - 1 - j)(n_{2} - 1 - i)} \right) | i \in \{0,\dots,n_{2} - 1\}, j \in \{0,\dots,n_{1} - 1\}.
		\end{displaymath}
	\end{definition}
	
	\begin{remark}
		For the square matrix $ \boldsymbol{M} = \left( m_{ij} \right) | i \in \{0,\dots,n - 1\} \ni j $ the next equation holds
		\begin{displaymath}
			\boldsymbol{M}^\tau = \boldsymbol{J}\boldsymbol{M}^T\boldsymbol{J},
		\end{displaymath}
		where $ \boldsymbol{J} $ is the matrix with ones on the diagonal $ i + j = n - 1 $ and noughts on other places.
	\end{remark}

	\begin{remark}
		Obviously, determinants of $ \boldsymbol{M} $ and $ \boldsymbol{M}^\tau $ are equal.
	\end{remark}
	
	\begin{proposition} \label{On Matrix Antitranspose}
		Let $ \boldsymbol{A} = \boldsymbol{R}\left[ i_{0},\dots,i_{u - 1}; j_{0},\dots,j_{v - 1} \right] $ be the submatrix of the matrix $ \boldsymbol{R} $ of dimension $ n \times n $, where $ u \leq n \geq v $. Then
		\begin{displaymath}
			\boldsymbol{A}^\tau = \boldsymbol{R}\left[ (n - 1 - j_{v - 1}),\dots,(n - 1 - j_{0}); (n - 1 - i_{u - 1}),\dots,(n - 1 - i_{0}) \right].
		\end{displaymath}
	\end{proposition}
	\begin{proof}
		\begin{displaymath}
			\boldsymbol{A} = \begin{pmatrix}
				a_{i_{0}j_{0}} & \cdots & a_{i_{0}j_{s}} & \cdots & a_{i_{0}j_{v-1}} \\
				\vdots & & & & \\
				a_{i_{r}j_{0}} & \cdots & a_{i_{r}j_{s}} & \cdots & a_{i_{r}j_{v-1}} \\ 
				\vdots & & & & \\
				a_{i_{u-1}j_{0}} & \cdots & a_{i_{u-1}j_{s}} & \cdots & a_{i_{u-1}j_{v-1}} \\ 
			\end{pmatrix} \Rightarrow \boldsymbol{A}^\tau \stackrel{def}{=} \begin{pmatrix}
				a_{i_{u-1}j_{v-1}} & \cdots & a_{i_{r}j_{v-1}} & \cdots & a_{i_{0}j_{v-1}} \\
				\vdots & & & & \\
				a_{i_{u-1}j_{s}} & \cdots & a_{i_{r}j_{s}} & \cdots & a_{i_{0}j_{s}} \\ 
				\vdots & & & & \\
				a_{i_{u-1}j_{0}} & \cdots & a_{i_{r}j_{0}} & \cdots & a_{i_{0}j_{0}} \\ 
			\end{pmatrix}
		\end{displaymath}
		Note, that $ a_{ij} = a_{(j - i + 1)_{\bmod{n}}} = a_{\left( (n - 1 - i) - (n - 1 - j) + 1 \right)_{\bmod{n}} } = a_{(n - 1 - j)(n - 1 - i)} $ for each element $ a_{ij} $ of $ \boldsymbol{R} $, therefore
		\begin{align*}
			\boldsymbol{A}^\tau \stackrel{def}{=}& \begin{pmatrix}
				a_{i_{u-1}j_{v-1}} & \cdots & a_{i_{r}j_{v-1}} & \cdots & a_{i_{0}j_{v-1}} \\
				\vdots & & & & \\
				a_{i_{u-1}j_{s}} & \cdots & a_{i_{r}j_{s}} & \cdots & a_{i_{0}j_{s}} \\ 
				\vdots & & & & \\
				a_{i_{u-1}j_{0}} & \cdots & a_{i_{r}j_{0}} & \cdots & a_{i_{0}j_{0}} \\ 
			\end{pmatrix} \\ 
			=& \begin{pmatrix}
				a_{(n - 1 - j_{v-1})(n - 1 - i_{u-1})} & \cdots & a_{(n - 1 - j_{v-1})(n - 1 - i_{r})} & \cdots & a_{(n - 1 - j_{v-1})(n - 1 - i_{0})} \\
				\vdots & & & & \\
				a_{(n - 1 - j_{s})(n - 1 - i_{u-1})} & \cdots & a_{(n - 1 - j_{s})(n - 1 - i_{r})} & \cdots & a_{(n - 1 - j_{s})(n - 1 - i_{0})} \\ 
				\vdots & & & & \\
				a_{(n - 1 - j_{0})(n - 1 - i_{u-1})} & \cdots & a_{(n - 1 - j_{0})(n - 1 - i_{r})} & \cdots & a_{(n - 1 - j_{0})(n - 1 - i_{0})} \\ 
			\end{pmatrix} \\
			=& \boldsymbol{R}\left[ (n - 1 - j_{v-1}),\dots,(n - 1 - j_{0}); (n - 1 - i_{u-1}),\dots,(n - 1 - i_{0}) \right].
		\end{align*}
		Furthermore, 
		\begin{displaymath}
			\left\lbrace \begin{aligned}
				i_{0} < \dots < &i_{r} < \dots < i_{u - 1}; \\
				j_{0} < \dots < &j_{s} < \dots < j_{v - 1},
			\end{aligned} \right.  \Longleftrightarrow \left\lbrace \begin{aligned}
				n - 1 - i_{u - 1} < \dots < &n - 1 - i_{r} < \dots < n - 1 - i_{0}; \\
				n - 1 - j_{v - 1} < \dots < &n - 1 - j_{s} < \dots < n - 1 - j_{0},
			\end{aligned} \right. 
		\end{displaymath}
		hence $ \boldsymbol{A}^\tau $ indeed is the submatrix of $ \boldsymbol{R} $.
	\end{proof}
	\begin{corollary}
		Under the conditions of theorem \ref{On Matrix Antitranspose}
		\begin{displaymath}
			\left\lbrace \boldsymbol{A} = \boldsymbol{A}^\tau \right\rbrace \Longleftrightarrow \left\lbrace \{ i_{0},\dots,i_{u - 1} \} = \{ n - 1 - j_{v - 1},\dots, n - 1 - j_{0} \} \right\rbrace. 
		\end{displaymath}
	\end{corollary}
	\begin{proof}
		To prove necessity suppose that $ \boldsymbol{A} = \boldsymbol{A}^\tau $. In this case $ u = v $, and $ a_{i_{r}j_{s}} = (\boldsymbol{A})_{rs} = \left( \boldsymbol{A}^\tau \right) _{rs} = a_{(n - 1 - j_{u - 1 - r})(n - 1 - i_{u - 1 - s})} $ for $ r \in \{0,\dots,u - 1\} $ and $ s \in \{0,\dots,u - 1\} $. As a corollary to remark \ref{On elements equality in L and R} the equality $ a_{i_{r}j_{s}} = a_{(n - 1 - j_{u - 1 - r})(n - 1 - i_{u - 1 - s})} $ of the elements of $ \boldsymbol{R} $ leads to the next system of equations
		\begin{displaymath}
			\left\lbrace \begin{aligned}
				i_{r} &= n - 1 - j_{u - 1 - r} + k \pmod n; \\ j_{s} &= n - 1 - i_{u - 1 - s} + k \pmod n,
			\end{aligned} \right. \Rightarrow \left\lbrace \begin{aligned}
				i_{r} &= n - 1 - j_{u - 1 - r}; \\ k &= 0,
			\end{aligned} \right.			
		\end{displaymath}
		hence $ \{ i_{0},\dots,i_{u - 1} \} = \{ n - 1 - j_{v - 1},\dots,n - 1 - j_{0} \} $.
		
		To prove sufficiency suppose that $ \{ i_{0},\dots,i_{u - 1} \} = \{ n - 1 - j_{v - 1},\dots,n - 1 - j_{0} \} $. Then
		\begin{align*}
			\boldsymbol{A} &= \boldsymbol{R}\left[ i_{0},\dots,i_{u - 1}; j_{0},\dots,j_{v - 1} \right] \\ 
			&= \boldsymbol{R}\left[ (n - 1 - j_{v - 1}),\dots,(n - 1 - j_{0}); (n - 1 - i_{u - 1}),\dots,(n - 1 - i_{0}) \right] \\
			&= \boldsymbol{A}^\tau.
		\end{align*}
	\end{proof}
	
\section{Construction of set of circulant matrix submatrices}
	
	We now describe a method to construct a minimal set $ \mathscr{M} $ of submatrices of a circulant matrix, which is sufficient to testify whether the matrix is an MDS one or not. Since an MDS matrix cannot have nought elements one may avoid special care of the $ 1 \times 1 $ submatrices by a convention. 
	
	In essence the method constructs for each $ u \in \{2, \dots, n - 2\} $ a set $ \mathcal{I}_{u} $ of tuples of row ordinals and a set $ \mathcal{J}_{u} $ of tuples of column ordinals such that their Cartesian product $ \mathcal{I}_{u} \times \mathcal{J}_{u} $ generates all the pairwise non-equivalent square submatrices of order $ u $.
	
	In order to construct $ \mathcal{I}_{u} $, recall corollary \ref{Corollary of equivalence necessary and sufficient condition} of theorem \ref{Equivalence necessary and sufficient condition} and remark \ref{Corollary of equivalence necessary and sufficient condition for R}: tuples $ \boldsymbol{d}' $ and $ \boldsymbol{d} $, corresponding to the equivalent submatrices of a circulant matrix, can be formed from each other by applying circular shifts. Therefore, when $ \boldsymbol{d}' $ and $ \boldsymbol{d} $ do not have such an interrelation the corresponding submatrices are non-equivalent. As a result, construction of $ \mathcal{I}_{u} $ may be reduced to the construction of a set $ \mathscr{D} $ of the tuples $ \boldsymbol{d} = (d_{1}, \dots, d_{u}) $, which cannot be obtained from each other by applying circular shifts, with $ d_{1} + \cdots + d_{u} = n $. Note, that these tuples represent different compositions of $ n $ into exactly $ u $ positive parts.

	\begin{definition}
		A composition of a natural $ m $ into exactly $ l $ parts is a tuple $ (d_{1}, \dots, d_{l}) $ of positive integers such that $ d_{1} + \cdots + d_{l} = m $. The tuples that differ in the ordering of their elements define different compositions.
	\end{definition}

	Having the set $ \mathscr{D} $ of compositions of $ n $ into $ u $ parts, which cannot be obtained from each other by applying circular shifts, one can construct $ \mathcal{I}_{u} $ the following way: 	
	\begin{displaymath}
		\mathcal{I}_{u} = \left\lbrace \{ 0, d_{1}, (d_{1} + d_{2}), \dots, (d_{1} + \cdots + d_{u - 1}) \}: (d_{1}, \dots, d_{u}) \in \mathscr{D} \right\rbrace.
	\end{displaymath}
	Thereby, $ \mathcal{I}_{u} $ comprises only those sets of rows of $ \boldsymbol{L} $ or $ \boldsymbol{R} $, having $ i_{0} = 0 $.
	
	Now we shall highlight the issues of constructing the set $ \mathcal{J}_{u} $ of tuples of column ordinals. Having the tuple $ \{i_{0}, \dots, i_{u - 1}\} \in \mathcal{I}_{u} $ fixed, two distinct tuples $ \{j_{0}, \dots, j_{u - 1}\} $ and $ \{j'_{0}, \dots, j'_{u - 1}\} $ render two non-equivalent matrices because these matrices share at least one element with a different number of occurrences in each of them. Furthermore, by construction of $ \mathcal{I}_{u} $, two distinct tuples $ \{i_{0}, \dots, i_{u - 1}\} \in \mathcal{I}_{u} $ and $ \{i'_{0}, \dots, i'_{u - 1}\} \in \mathcal{I}_{u} $  render two non-equivalent matrices for each $ \{j_{0}, \dots, j_{u - 1}\} $ and $ \{j'_{0}, \dots, j'_{u - 1}\} $ due to corollary \ref{Corollary of equivalence necessary and sufficient condition}. Thereby, $ \mathcal{J}_{u} $ must comprise all the tuples possible in order that Cartesian product $ \mathcal{I}_{u} \times \mathcal{J}_{u} $ generates all the pairwise non-equivalent $ u \times u $ submatrices. In other words $ \mathcal{J}_{u} = \mathcal{P}_{u}(\{0,\dots,n - 1\}) $, where $ \mathcal{P}_{u}(\{0,\dots,n - 1\}) $ denotes the set of subsets of $ \{0,\dots,n - 1\} $ of cardinality $ u $.
	
	Apart from equivalent submatrices, there are transposes (in case of matrix $ \boldsymbol{L} $) or anti-transposes (in case of matrix $ \boldsymbol{R} $), which should be taken into account when constructing minimal $ \mathscr{M} $ set since their determinants equal to the determinant of the original submatrix. Propositions \ref{On Matrix Transpose} states that for a submatrix $ \boldsymbol{A}_{L} $ formed from rows $ \{i_{0}, \dots, i_{u - 1}\} $ and columns $ \{j_{0}, \dots, j_{u - 1} \} $ of $ \boldsymbol{L} $, the transpose $ \left( \boldsymbol{A}_{L} \right)^{T} $ is formed from rows $ \{j_{0}, \dots, j_{u - 1} \} $ and columns $ \{i_{0}, \dots, i_{u - 1}\} $ of $ \boldsymbol{L} $. Similarly, proposition \ref{On Matrix Antitranspose} states that for a submatrix $ \boldsymbol{A}_{R} $ formed from rows $ \{i_{0}, \dots, i_{u - 1}\} $ and columns $ \{j_{0}, \dots, j_{u - 1} \} $ of $ \boldsymbol{R} $, the anti-transpose $ \left( \boldsymbol{A}_{R} \right)^{\tau} $ is formed from rows $ \{(n - 1 - j_{v - 1}),\dots,(n - 1 - j_{0})\} $ and columns $ \{(n - 1 - i_{u - 1}),\dots,(n - 1 - i_{0})\} $ of the same matrix. It is noteworthy that two distinct tuples $ \boldsymbol{i} = \{i_{0}, \dots, i_{u - 1}\} $ and $ \boldsymbol{i}' = \{i'_{0}, \dots, i'_{u - 1}\} $ in $ \mathcal{I}_{u} $ reveal the following interrelation:
	\begin{multline*}
		\boldsymbol{L}\left[ i'_{0}, \dots, i'_{u - 1};  \pi\left( (i_{0} + k)_{\bmod n}\right) , \dots, \pi \left( (i_{u - 1} + k)_{\bmod n} \right) \right] \\
		= \left( \boldsymbol{L}\left[ \pi\left( (i_{0} + k)_{\bmod n}\right) , \dots, \pi \left( (i_{u - 1} + k)_{\bmod n} \right); i'_{0}, \dots, i'_{u - 1} \right] \right)^{T} \\
		\sim \left( \boldsymbol{L}\left[ i_{0}, \dots, i_{u - 1};  \pi'\left( (i'_{0} - k)_{\bmod n}\right) , \dots, \pi' \left( (i'_{u - 1} - k)_{\bmod n} \right) \right] \right)^{T}.
	\end{multline*}
	Here $ \pi $ is a permutation of $ \{(i_{0} + k)_{\bmod n}, \dots,  (i_{u - 1} + k)_{\bmod n}\} $ such that $ \pi\left( (i_{0} + k)_{\bmod n}\right) < \cdots < \pi\left( (i_{u - 1} + k)_{\bmod n}\right) $, while $ \pi' $ is a permutation of $ \{(i'_{0} - k)_{\bmod n}, \dots,  (i'_{u - 1} - k)_{\bmod n}\} $ such that $ \pi'\left( (i'_{0} - k)_{\bmod n}\right) < \cdots < \pi'\left( (i'_{u - 1} - k)_{\bmod n}\right) $, а $ k \in \{0, \dots, n - 1\} $.
	
	By analogy the next interrelation holds:
	\begin{multline*}
		\boldsymbol{R}\left[ i'_{0}, \dots, i'_{u - 1};  \pi\left( (i_{0} + k)_{\bmod n}\right) , \dots, \pi \left( (i_{u - 1} + k)_{\bmod n} \right) \right] \\
		= \left( \boldsymbol{R}\left[ \pi\left( (n - 1 - i_{u - 1} + k)_{\bmod n}\right) , \dots, \pi \left( (n - 1 - i_{0} + k)_{\bmod n} \right)\right. \right.; \\ 
		\left. \left. (n - 1 - i'_{u - 1}),\dots,(n - 1 - i'_{0}) \right] \right)^{\tau} \\
		\sim \left( \boldsymbol{R}\left[ i_{0}, \dots, i_{u - 1};  \pi'\left( (i_{0} + i_{u - 1} - i'_{u - 1} - k)_{\bmod n}\right) , \dots, \pi' \left( (i_{0} + i_{u - 1} - i'_{0} - k)_{\bmod n} \right) \right] \right)^{\tau}.
	\end{multline*}
	Here $ \pi $ is a permutation of $ \{(i_{0} + k)_{\bmod n}, \dots,  (i_{u - 1} + k)_{\bmod n}\} $ such that $ \pi\left( (i_{0} + k)_{\bmod n}\right) < \cdots < \pi\left( (i_{u - 1} + k)_{\bmod n}\right) $, while $ \pi' $ is a permutation of $ \{(i_{0} + i_{u - 1} - i'_{u - 1} - k)_{\bmod n}, \dots,  (i_{0} + i_{u - 1} - i'_{0} - k)_{\bmod n}\} $ such that
	\begin{displaymath}
		\pi'\left( (i_{0} + i_{u - 1} - i'_{u - 1} - k)_{\bmod n}\right) < \cdots < \pi'\left( (i_{0} + i_{u - 1} - i'_{0} - k)_{\bmod n}\right),
	\end{displaymath}
	and $ k \in \{0, \dots, n - 1\} $.
	
	Thus, on the basis of the above observations, one can conclude that to avoid the surplus inclusion of transposes or anti-transposes which are taken into account when processing the preceding tuple $ \boldsymbol{i} = \{i_{0}, \dots, i_{u - 1}\}$, it is needed to exclude for every $ k \in \{0, \dots, n - 1\} $ the tuples $ \left\lbrace \pi\left( (i_{0} + k)_{\bmod n}\right), \dots, \pi \left( (i_{u - 1} + k)_{\bmod n} \right) \right\rbrace $ from $ \mathcal{P}_{u}(\{0,\dots,n - 1\}) $  when processing the next tuple $ \boldsymbol{i}' $ within construction of the set $ \mathcal{J}_{u} $.
	
	We now give a more formal description of the algorithm that constructs for matrices $ \boldsymbol{L} $ or $ \boldsymbol{R} $ the set $ \mathscr{M} $ of all non-equivalent submatrices without equivalence classes of transposes or anti-transposes.

	\begin{algorithm}[H]
		\caption{Construction of a set of all non-equivalent submatrices of a circulant matrix without equivalence classes of transposes or anti-transposes}
		\label{Matrix minimal matrices set construction}
		\begin{algorithmic}[1]
			\Require A circulant matrix $ \boldsymbol{M} = \boldsymbol{L} $ or $ \boldsymbol{M} = \boldsymbol{R} $ of the order $ n > 2 $. 
			\Ensure A set $ \mathscr{M} $ of all non-equivalent submatrices without equivalence classes of transposes or anti-transposes
			\State $ \mathscr{M} \gets \emptyset  $
			\ForAll{$ u \in \{2, \dots, n - 1\} $}
				\State \parbox[t]{\dimexpr\linewidth-\algorithmicindent} {Construct a set $ \mathscr{D}_{\prec} $ of compositions of $ n $ into $ u $ positive parts equipped with a linear order $ \prec $ fixed for all $ u $:}
					\begin{displaymath}
						\mathscr{D}_{\prec} \gets \left\lbrace (d_{1}, \dots, d_{u}): d_{1} + \cdots + d_{u} = m, d_{l} > 0, l \in \{0, \dots, n - 1\} \right\rbrace
					\end{displaymath}
%				\State $ \mathscr{D}_{i} \gets \left\lbrace (d_{1}, \dots, d_{u}): d_{1} + \cdots + d_{u} = m \right\rbrace $ \Comment{Построить множество $ \mathscr{D}_{i} $ композиций числа $ n $ длины $ u $ без нулевых координат}
				\State \parbox[t]{\dimexpr\linewidth-\algorithmicindent} {For each $ \boldsymbol{d} \in \mathscr{D}_{\prec} $ find and eliminate from $ \mathscr{D}_{\prec} $ circular shifts of $ \boldsymbol{d} $ distinct from $ \boldsymbol{d} $}
				\State \parbox[t]{\dimexpr\linewidth-\algorithmicindent} {Construct a set $ \mathcal{I}_{u} $ of tuples of row ordinals in $ \boldsymbol{M} $ of non-equivalent $ u \times u $ submatrices:
					\begin{displaymath}
						\mathcal{I}_{u} \gets \left\lbrace \{ 1, (1 + d_{1}), (1 + d_{1} + d_{2}), \dots, (1 + d_{1} + \cdots + d_{u - 1}) \}: (d_{1}, \dots, d_{u}) \in \mathscr{D}_{\prec} \right\rbrace
					\end{displaymath}
				}
%				\State $ \mathscr{I} \gets \left\lbrace \{ 1, (1 + d_{1}), (1 + d_{1} + d_{2}), \dots, (1 + d_{1} + \cdots + d_{u - 1}) \}| (d_{1}, \dots, d_{u}) \in \mathscr{D}_{i} \right\rbrace $ 
%				\Comment{ Построить множество $ \mathscr{I} $ номеров строк неэквивалентных подматриц размера $ u \times u $}
				\State \parbox[t]{\dimexpr\linewidth-\algorithmicindent} {Construct a set $ \mathcal{J}_{u} $ of tuples of column ordinals in $ \boldsymbol{M} $ of non-equivalent $ u \times u $ submatrices:
					\begin{displaymath}
						\mathcal{J}_{u} \gets \mathcal{P}_{u}(\{0,\dots,n - 1\}),
					\end{displaymath}
					where $ \mathcal{P}_{u}(\{0,\dots,n - 1\}) $ be the set of the subsets of $ \{0,\dots,n - 1\} $ of cardinality $ u $}
%				\State $ \mathscr{J} \gets \mathcal{P}_{u}(\{0,\dots,n - 1\}) $
%				\Comment{$ \mathcal{P}_{u}(\{0,\dots,n - 1\}) $ --- множество всех $ u $-элементных подмножеств множества $ \{0,\dots,n - 1\} $}
%				\Comment{Построить множество $ \mathscr{J} $ номеров столбцов неэквивалентных подматриц размера $ u \times u $}
				\State $ \mathcal{J}'_{u} \gets \emptyset $
				\ForAll {$ \boldsymbol{i} \in \mathcal{I}_{u} $}
				\State $ \mathcal{J}_{u} \gets \mathcal{J}_{u} \setminus \mathcal{J}'_{u} $
					\ForAll {$ \boldsymbol{j} \in \mathcal{J}_{u} $}
						\State $ \mathscr{M} \gets \mathscr{M} \cup \boldsymbol{M}\left[i_{0}, \dots, i_{u - 1}; j_{0}, \dots, j_{u - 1}\right] $: $ \{i_{0}, \dots, i_{u - 1}\} = \boldsymbol{i}, \{j_{0}, \dots, j_{u - 1}\} = \boldsymbol{j} $
					\EndFor
					\State \parbox[t]{\dimexpr\linewidth-\algorithmicindent-\algorithmicindent} {Add the set of tuples of row ordinals of non-equivalent submatrices to $ \mathcal{J}'_{u} $:
						\begin{displaymath}
							\mathcal{J}'_{u} \gets \mathcal{J}'_{u} \cup \left\lbrace \left\lbrace \pi\left( (i_{0} + k)_{\bmod n}\right) , \dots, \pi \left( (i_{u - 1} + k)_{\bmod n} \right) \right\rbrace  | k \in \{0, \dots, n - 1\} \right\rbrace,
						\end{displaymath}
						where $ \pi $ be a permutation of $ \{(i_{0} + k)_{\bmod n}, \dots,  (i_{u - 1} + k)_{\bmod n}\} $ such that $ \pi\left( (i_{0} + k)_{\bmod n}\right) < \cdots < \pi\left( (i_{u - 1} + k)_{\bmod n}\right), k \in \{0, \dots, n - 1\} $}
%					\State $ \mathscr{J}' \gets \mathscr{J}' \cup \left\lbrace \left\lbrace \pi\left( (i_{0} + k)_{\bmod n}\right) , \dots, \pi \left( (i_{u - 1} + k)_{\bmod n} \right) \right\rbrace  | k \in \{0, \dots, n - 1\} \right\rbrace $, где $ \pi $ --- перестановка элементов множества $ \{(i_{0} + k)_{\bmod n}, \dots,  (i_{u - 1} + k)_{\bmod n}\} $ такая, что $ \pi\left( (i_{0} + k)_{\bmod n}\right) < \cdots < \pi\left( (i_{u - 1} + k)_{\bmod n}\right) $
				\EndFor
			\EndFor
			\State $ \mathscr{M} \gets \mathscr{M} \cup \boldsymbol{M}$
			\State \Return $ \mathscr{M} $	
		\end{algorithmic}
	\end{algorithm}
%\newpage

We shall further estimate a power of the set $ \mathscr{M} $ that algorithm \ref{Matrix minimal matrices set construction} constructs and compare it to the overall quantity of square submatrices. There are two cases possible, namely $ n $ is a prime number or a composite one.

Consider first a case $ n $ be a prime. By remark \ref{Remark on prime order matrix submatrices quantity}, every generator submatrix of a circulant matrix $ \boldsymbol{M} $ of the prime order forms an equivalence class with $ n $ elements. Consequently, the number of equivalence classes is
\begin{displaymath}
	1 + \frac{1}{n} \sum_{u = 1}^{n - 1}\binom{n}{u}^2,
\end{displaymath}
where the term "1" \;\;stands for the only square submatrix of order $ n $, and the binomial coefficient squared is equal to the number of the square submatrices of order $ u $. Besides, the method discussing eliminates transposes or anti-transposes, which again form equivalence classes with $ n $ elements. At the same time there exist symmetric submatrices $ \boldsymbol{A}_{M} $ of $ \boldsymbol{M} $, i.e. $ \boldsymbol{A}_{M} = \boldsymbol{A}^{T}_{M} $ and $ \boldsymbol{A}^{T}_{M} $ is a submatrix of $ \boldsymbol{M} $ itself. The quantity of symmetric submatrices of order $ u $ equals the number of subsets of cardinality $ u $ of an n element set. Then, taking transposes of the submatrices and symmetric submatrices into account, 
\begin{displaymath}
	\left| \mathscr{M} \right| = 1 + \sum_{u = 1}^{n - 1}\binom{n}{u} + \frac{1}{2 n} \sum_{u = 1}^{n - 1}\left( \binom{n}{u}^2 - n \binom{n}{u} \right),
\end{displaymath}
where the second term denotes the number of symmetric submatrices which generate distinct equivalence classes, and the third term determine the number of equivalence classes that do not comprise symmetric submatrices.

Consider now a matrix $ \boldsymbol{M} $ of a composite order $ n $. By virtue of proposition \ref{k-values for A_L} and remark \ref{Remark on A_R to k-values for A_L}, an equivalence class generated by some square matrix of order $ u $ comprises $ k \in \left\lbrace qn / \gcd(n, u) | q \in \left\lbrace 1, \dots, \gcd(n, u) \right\rbrace \right\rbrace $ elements. It can be seen that for coprime $ u $ and $ n $ the equivalence class comprises exactly $ n $ elements, but $ n / \gcd(n, u) \leq k \leq n $ when $ \gcd(n, u) > 1 $. Thus, 
\begin{displaymath}
	\left| \mathscr{M} \right| \leq 1 + \sum_{u=1}^{n-1}\binom{n}{u} + \frac{1}{2 n} \sum_{u=1}^{n-1}\left( \gcd(n,u)\binom{n}{u}^2 - n\binom{n}{u} \right).
\end{displaymath}

Thereby, in general case
\begin{displaymath}
	\frac{1}{2 n} \sum_{u = 1}^{n - 1} \binom{n}{u}^2 \leq \left| \mathscr{M} \right| + \frac{1}{2} \sum_{u = 1}^{n - 1}\binom{n}{u} - 1 \leq \frac{1}{2 n} \sum_{u=1}^{n-1} \gcd(n,u)\binom{n}{u}^2.
\end{displaymath}
	
It turns out that the determinant of each submatrix in $ \mathscr{M} $ can be expanded by determinants of the submatrices which are again in $ \mathscr{M} $. The next proposition proves this fact.
\begin{proposition}
	The determinant of submatrix of order $ u $ in the set $ \mathscr{M} $ constructed by algorithm \ref{Matrix minimal matrices set construction} can be expanded by the square submatrices of order $ u - 1 $, lying in $ \mathscr{M} $ themselves. 
\end{proposition}
\begin{proof}
	One can prove that the last row determinant expansion can be performed utilizing determinants of only the submatrices in the set $ \mathscr{M} $, which is constructed by algorithm \ref{Matrix minimal matrices set construction}.  
	
	Let $ \boldsymbol{M}\left[i_{0}, \dots, i_{u - 1}; j_{0}, \dots, j_{u - 1}\right] \in \mathscr{M} $ be a submatrix of order $ u $, and observe a submatrix $ \boldsymbol{M}\left[i_{0}, \dots, i_{u - 2}; j_{0}, \dots, j_{s-1}, j_{s+1}, \dots, j_{u - 1}\right] $ formed by deleting the row $ i_{u-1} $ and the column $ j_{s}, j_{0} \leq j_{s} \leq j_{u-1} $. Note, that the submatrices formed from the rows $ (i_{0}, \dots, i_{u - 2}) $ are in $ \mathscr{M} $. Otherwise, the submatrices formed from the rows $ (i'_{0}, \dots, i'_{u - 2}) \prec (i_{0}, \dots, i_{u - 2}) $, for which $ i'_{1} - i'_{0} = i_{1} - i_{0}, \dots, i'_{u-2} - i'_{u-3} = i_{u-2} - i_{u-3} $, would have belonged to $ \mathscr{M} $, while $ \boldsymbol{M}\left[i_{0}, \dots, i_{u - 1}; j_{0}, \dots, j_{u - 1}\right] $ would have escaped $ \mathscr{M} $ (the submatrix formed from the rows $ (i'_{0}, \dots, i'_{u - 2}, i'_{u-2} + i_{u-1} - i_{u-2}) $ would have been instead in $ \mathscr{M} $). 
	
	Moreover, $ \boldsymbol{M}\left[i_{0}, \dots, i_{u - 2}; j_{0}, \dots, j_{s-1}, j_{s+1}, \dots, j_{u - 1}\right] \in \mathscr{M} $, otherwise a tuple $ (i'_{0}, \dots, i'_{u - 2}) $ such that $ (j_{0}, \dots, j_{s-1}, j_{s+1}, \dots, j_{u - 1}) = \left( \pi\left( (i'_{0} + k)_{\bmod n}\right), \dots, \pi\left( (i'_{u - 2} + k)_{\bmod n}\right) \right) $ for some $ k \in \{0, \dots, n - 1\} $ would have existed. Here $ \pi $\; is a permutation that arranges its arguments in ascending order. Note that $ (i'_{0}, \dots, i'_{u - 2}) \prec (i_{0}, \dots, i_{u - 2}) $, while $ (i'_{0}, \dots, i'_{s-1}, i'_{s-1} + j_{s} - j_{s-1}, i'_{s+1}, \dots, i'_{u - 2}) \prec (i_{0}, \dots, i_{u - 1}) $, since $ i'_{0} \leq i_{0}, \dots, i'_{u-2} < i_{u-2} < i_{u-1}$. Thereby, during construction of the set $ \mathcal{J}_{u} $ all the elements
	\begin{multline*}
		\left( \pi\left( (i'_{0} + k)_{\bmod n}\right), \dots, \pi\left( (i'_{s-1} + k)_{\bmod n} \right), \pi\left( (i'_{s-1} + j_{s}- j_{s-1} + k)_{\bmod n} \right),\right. \\ \left.  \pi\left( (i'_{s+1} + k)_{\bmod n}\right), \dots, \pi\left( (i_{u - 2} + k)_{\bmod n}\right) \right)
	\end{multline*}
	would have been deleted from  $ \mathcal{P}_{u}(\{0,\dots,n - 1\}) $, and the tuple $ (j_{0}, \dots, j_{u - 1}) $ in particular. This results in $ \boldsymbol{M}\left[i_{0}, \dots, i_{u - 1}; j_{0}, \dots, j_{u - 1}\right] \notin \mathscr{M} $, leading to a contradiction.
	
	Thus, for each submatrix of order $ u $ in $ \mathscr{M} $ a submatrix of order $ u - 1 $, formed by deleting the last row and some column, lies in $ \mathscr{M} $ itself.
\end{proof}

We now present practical results regarding the cardinality of $ \mathscr{M} $ for several circulant matrices of order $ n \leq 16 $. The next table \ref{Experiments of M power rate} gives for different values of $ n $ the overall quantity $ Q $ of square submatrices, the cardinality of the corresponding set $ \mathscr{M} $ and a rate $ Q / |\mathscr{M}| $.

\begin{table}[H]
	\caption{The ratio of an overall quantity $ Q $ of square submatrices of a circulant matrix of order $ n $ to the cardinality of $ \mathscr{M} $}
	\begin{center} \label{Experiments of M power rate}
		\begin{tabular}[c]{|c||c|c|c|c|}
			\hline
			$ n $ & 7 & 8 & 12 & 16 \\
			\hline
			\hline
			$ Q $ & 3431 & 12869 & 2 704 155 & 601 080 389 \\
			\hline
			$ |\mathscr{M}| $ & 309 & 941 & 115 157 & 18 838 305 \\
			\hline
			$ Q / |\mathscr{M}| $ & 11.1 & 13.7 & 23.5 & 31.9 \\
			\hline
		\end{tabular}
	\end{center}
\end{table}

\section{Conclusion}

In this paper we focused on circulant matrices and proposed several techniques which allow reduction of computational efforts into circulant MDS matrix verification. This was possible due to the fact that circulant matrices have sets of submatrices with up-to-a-sign-equal determinants. This fact allows avoiding the necessity to compute all the submatrices determinants in order to perform MDS verification. Taking equivalent submatrices, their transposes or anti-transposes into account, we were able to construct approximately $ 2n $ times smaller set of submatrices needing for MDS verification, where $ n $ is the order of a circulant matrix. We also gave the explicit description of the algorithm that constructs the set in question.
	
\printbibliography

\end{document}